\input ./style/arxiv-general.cfg
\documentclass[aos,MSNbibl,nameyear,dvips]{arximspdf}
\makeatletter
   \@ifpackageloaded{graphicx}{}{\usepackage{graphicx}}
\makeatother

\doi{10.1214/15-AOS1384}
\volume{44}
\issue{2}
\pubyear{2016}
\firstpage{713}
\lastpage{742}
\docsubty{FLA}

\makeatletter

\newcommand{\rrvert}{\vert}
\newcommand{\llvert}{\vert}
\def\cal{\mathcal}

\newproclaim{definition}{Definition}

\newproclaim{example}{Example}
\newtheorem{theorem}{Theorem}
\newtheorem{corollary}[theorem]{Corollary}
\newtheorem{lemma}[theorem]{Lemma}

\newcommand{\openr}{\mathbb{R}}
\newcommand{\pr}{\operatorname{Pr}_0}
\newcommand{\Var}{\operatorname{Var}}
\makeatother

\begin{document}
\begin{frontmatter}

\title{Statistical inference for the mean outcome under a possibly
non-unique optimal treatment strategy}
\runtitle{Non-unique optimal treatment strategies}

\begin{aug}
\author[A]{\fnms{Alexander R.}~\snm{Luedtke}\corref{}\thanksref{T1,T2}\ead
[label=e1]{aluedtke@berkeley.edu}}
\and
\author[A]{\fnms{Mark~J.}~\snm{van der Laan}\thanksref{T1}\ead
[label=e2]{laan@berkeley.edu}}
\runauthor{A.~R. Luedtke and M.~J. van der Laan}
\thankstext{T1}{Supported by NIH Grant R01 AI074345-06.}
\thankstext{T2}{Supported by the Department of Defense (DoD) through
the National Defense Science \& Engineering Graduate Fellowship (NDSEG)
Program.}
\affiliation{University of California, Berkeley}
\address[A]{Division of Biostatistics\\
University of California, Berkeley\\
101 Haviland Hall\\
Berkeley, California 94720\\
USA\\
\printead{e1}\\
\phantom{E-mail:\ }\printead*{e2}}
\end{aug}

%
\received{\smonth{12} \syear{2014}}
%
\revised{\smonth{9} \syear{2015}}


\begin{abstract}
We consider challenges that arise in the estimation of the mean outcome
under an optimal individualized treatment strategy defined as the
treatment rule that maximizes the population mean outcome, where the
candidate treatment rules are restricted to depend on baseline
covariates. We prove a necessary and sufficient condition for the
pathwise differentiability of the optimal value, a key condition needed
to develop a regular and asymptotically linear (RAL) estimator of the
optimal value. The stated condition is slightly more general than the
previous condition implied in the literature. We then describe an
approach to obtain root-$n$ rate confidence intervals for the optimal
value even when the parameter is not pathwise differentiable. We
provide conditions under which our estimator is RAL and asymptotically
efficient when the mean outcome is pathwise differentiable. We also
outline an extension of our approach to a multiple time point problem.
All of our results are supported by simulations.
\end{abstract}

%

\end{frontmatter}

\section{Introduction}\label{sec1}
There has been much recent work in estimating optimal treatment regimes
(TRs) from a random sample. A TR is an individualized treatment
strategy in which treatment decisions for a patient can be based on
their measured covariates. Doctors generally make decisions this way,
and thus it is natural to want to learn about the best strategy. The
value of a TR is defined as the population counterfactual mean outcome
if the TR were implemented in the population. The optimal TR is the TR
with the maximal value, and the value at the optimal TR is the optimal
value. In a single time point setting, the optimal TR can be defined as
the sign of the ``blip function,'' defined as the additive effect of a
blip in treatment on a counterfactual outcome, conditional on baseline
covariates [\citet{Robins04}]. In a multiple time point setting,
treatment strategies are called dynamic TRs (DTRs). For a general
overview of recent work on optimal (D)TRs, see \citet{ChakrabortyMoodie2013}.

Suppose one wishes to know the impact of implementing an optimal TR in
the population, that is, one wishes to know the optimal value. Before
estimating the optimal value, one typically estimates the optimal rule.
Recently, researchers have suggested applying machine learning
algorithms to estimate the optimal rules from large classes which
cannot be described by a finite dimensional parameter [see,
e.g., \citet{Zhangetal12,Zhaoetal12,LuedtkevanderLaan2014}].

Inference for the optimal value has been shown to be difficult at
exceptional laws, that is, probability distributions where there exists
a strata of the baseline covariates that occurs with positive
probability and for which treatment is neither beneficial nor harmful
[\citet{Robins04,RobinsRotnitzky14}]. \citet{Zhangetal12b} considered
inference for the optimal value in restricted classes in which the TRs
are indexed by a finite-dimensional vector. At non-exceptional laws,
they outlined an argument showing that their estimator is (up to a
negligible term) equal to the estimator that estimates the value of the
\textit{known} optimal TR under regularity conditions. The implication
is that one can estimate the optimal value and then use the usual
sandwich technique to estimate the standard error and develop Wald-type
confidence intervals (CIs). \citet{vanderLaanLuedtke2014} and \citet
{vanderLaanLuedtke14b} developed inference for the optimal value when
the DTR belongs to an unrestricted class. \citet{vanderLaanLuedtke14b}
provide a proof that the efficient influence curve for the parameter
which treats the optimal rule as known is equal to the efficient
influence curve of the optimal value at non-exceptional laws. One of
the contributions of the current paper is to present a slightly more
precise statement of the condition for the pathwise differentiability
of the mean outcome under the optimal rule. We will show that this
condition is necessary and sufficient.

However, restricting inference to non-exceptional laws is limiting as
there is often no treatment effect for people in some strata of
baseline covariates. \citet{Chakrabortyetal2014} propose using the
$m$-out-of-$n$ bootstrap to obtain inference for the value of an
estimated DTR. With an inverse probability weighted (IPW) estimator
this yields valid inference when the treatment mechanism is known or is
estimated according to a correctly specified parametric model. They
also discuss an extension to an double robust estimator. The
$m$-out-of-$n$ bootstrap draws samples of size $m$ patients from the
data set of size $n$. In non-regular problems, this method yields valid
inference if $m,n\rightarrow\infty$ and $m=o(n)$. The CIs for the value
of an estimated regime shrink at a root-$m$ (not root-$n$) rate. In
addition to yielding wide CIs, this approach has the drawback of
requiring a choice of the important tuning parameter $m$, which
balances a trade-off between coverage and efficiency. \citeauthor
{Chakrabortyetal2014} propose using a double bootstrap to select this
tuning parameter.


\citet{Goldbergetal2014} instead consider truncating the criteria to be
optimized, that is, the value under a given rule, so that only
individuals with a clinically meaningful treatment effect contribute to
the objective function. These authors then propose proceeding with
inference for the truncated value at the optimal DTR. For a fixed
truncation level, the estimated truncated optimal value minus the true
truncated optimal value, multiplied by root-$n$, converges to a normal
limiting distribution. \citet{Laberetal2014b} propose instead replacing
the indicator used to define the value of a TR with a differentiable
function. They discuss situations in which the estimator minus the
smoothed value of the estimated TR, multiplied by root-$n$, would have
a reasonable limit distribution.

In this work, we develop root-$n$ rate inference for the optimal value
under reasonable conditions. Our approach avoids any sort of
truncation, and does not require that the estimate of the optimal rule
converge to a fixed quantity as the sample size grows. We show that our
estimator minus the truth, properly standardized, converges to a
standard normal limiting distribution. This allows for the
straightforward construction of asymptotically valid CIs for the
optimal value. Neither the estimator nor the inference rely on a
complicated tuning parameter. We give conditions under which our
estimator is asymptotically efficient among all regular and
asymptotically linear (RAL) estimators when the optimal value parameter
is pathwise differentiable, similar to those we presented in \citet
{vanderLaanLuedtke2014}. However, they do not require that one knows
that the optimal value parameter is pathwise differentiable from the
outset. Implementing the procedure only requires a minor modification
to a typical one-step estimator.

We believe the value of the unknown optimal rule is an interesting
target of inference because the treatment strategy learned from the
given data set is likely to be improved upon as clinicians gain more
knowledge, with the treatment strategy given in the population
eventually approximating the optimal rule. Additionally, the optimal
rule represents an upper bound on what can be hoped for when a
treatment is introduced. Nonetheless, as we and others have argued in
the references above, the value of the estimated rule is also an
interesting target of inference [\citet
{Chakrabortyetal2014,Laberetal2014b},
\citeauthor{vanderLaanLuedtke14b}
(\citeyear{vanderLaanLuedtke14b,vanderLaanLuedtke2014})].
Thus, although our focus is on estimating the optimal value, we also
give conditions under which our CI provides proper coverage for the
data adaptive parameter which gives the value of the rule estimated
from the entire data set. 


\subsection*{Organization of article}
Section~\ref{sec:formulation} formulates the statistical problem of
interest. Section~\ref{sec:pd} gives necessary and sufficient conditions
for the pathwise differentiability of the optimal value. Section~\ref
{sec:challenge} outlines the challenge of obtaining inference at
exceptional laws and gives a thought experiment that motivates our
procedure for estimating the optimal value. Section~\ref{sec:martcis}
presents an estimator for the optimal value. This estimator represents
a slight modification to a recently presented online one-step estimator
for pathwise differentiable parameters. Section~\ref{sec:compeff} discusses
computationally efficient implementations of our proposed procedure.
Section~\ref{sec:discconds} discusses each condition of the key result
presented in Section~\ref{sec:martcis}. Section~\ref{sec:sims}
describes our
simulations. Section~\ref{sec:simresults} gives our simulation results.
Section~\ref{sec:disc} closes with a summary and some directions for
future work.

All proofs can be found in Supplementary Appendix A [\citet
{LuedtkevanderLaan2014b}]. We outline an extension of our proposed
procedure to the multiple time point setting in Supplementary Appendix
B. Additional figures appear in Supplementary Appendix C.

\section{Problem formulation} \label{sec:formulation}
Let $O=(W,A,Y)\sim P_0\in{\cal M}$, where $W$ represents a vector of
covariates, $A$ a binary intervention, and $Y$ a real-valued outcome.
The model for $P_0$ is nonparametric. We observe an independent and
identically distributed (i.i.d.) sample $O_1,\ldots,O_n$ from $P_0$. Let
$\mathcal{W}$ denote the range of $W$. For a distribution $P$, define
the treatment mechanism $g(P)(A|W)\triangleq\operatorname{Pr}_P(A|W)$.
We will refer to $g(P_0)$ as $g_0$ and $g(P)$ as $g$. For a function
$f$, we will use $E_P[f(O)]$ to denote $\int f(o)\,dP(o)$. We will also
use $E_{0}[f(O)]$ to denote $E_{P_0}[f(O)]$ and $\pr$ to denote the
$P_0$ probability of an event.
Let $\Psi:{\cal M}\rightarrow\openr$ be defined by
\[
\Psi(P)\triangleq E_P E_P\bigl[Y|A=d(P) (W),W\bigr],
\]
where $d(P)\triangleq\operatorname{argmax}_{d} E_PE_P(Y|A=d(W),W)$ is an
optimal treatment rule under $P$. We will resolve the ambiguity in the
definition of $d$ when the $\operatorname{argmax}$ is not unique later in
this section. Throughout we assume that $\pr(0<g_0(1|W)<1)$ so that
$\Psi(P_0)$ is well defined. Under causal assumptions, $\Psi(P)$ is
equal to the counterfactual mean outcome if, possibly contrary to fact,
the rule $d(P)$ were implemented in the population. We can also
identify $d(P)$ with a causally optimal rule under those same
assumptions. We refer the reader to \citet{vanderLaanLuedtke2014} for a
more precise formulation of such a treatment strategy. As the focus of
this work is statistical, all of the results will hold when estimating
the parameter $\Psi(P_0)$ whether or not the causal assumptions needed
for identifiability hold. Let
\begin{eqnarray*}
\bar{Q}(P) (A,W)&\triangleq& E_P[Y|A,W],
\\
\bar{Q}_b(P) (W)&\triangleq&\bar{Q}(P) (1,W)-\bar{Q}(P) (0,W).
\end{eqnarray*}
We will refer to $\bar{Q}_b(P)$ the blip function for the distribution
$P$. We will denote to the above quantities applied to $P_0$ as $\bar
{Q}_0$ and $\bar{Q}_{b,0}$, respectively. We will often omit the
reliance on $P$ altogether when there is only one distribution $P$
under consideration: $\bar{Q}(A,W)$ and $\bar{Q}_b(W)$. We also define
$\Psi_d(P)= E_P \bar{Q}(d(P)(W),W)$. Consider the efficient influence
curve of $\Psi_d$ at $P$:
\[
D(d,P) (O)=\frac{I(A=d(W))}{g(A|W)}\bigl(Y-\bar{Q}(A,W)\bigr)+\bar {Q}\bigl(d(W),W
\bigr)-\Psi_d(P).
\]
Let $B(P)\triangleq\{w: \bar{Q}_b(w)=0\}$. We will refer to $B(P_0)$
as $B_0$. An exceptional law is defined as a distribution $P$ for which
$\operatorname{Pr}_P(W\in B(P))>0$ [\citet{Robins04}]. We note that the
ambiguity in the definition of $d(P)$ occurs precisely on the set
$B(P)$. In particular, $d(P)$ must almost surely agree with some rule
in the class
%
\begin{equation}
\bigl\{w\mapsto I\bigl(\bar{Q}_{b}(w)>0\bigr)I\bigl(w\notin B(P)
\bigr) + b(w)I\bigl(w\in B(P)\bigr) : b \bigr\}, \label{eq:optclass}
\end{equation}
where $b : \mathcal{W}\rightarrow\{0,1\}$ is some function. Consider
now the following uniquely defined optimal rule:
\[
d^*(P) (W)\triangleq I\bigl(\bar{Q}_b(W)>0\bigr).
\]
We will let $d_0^*=d^*(P_0)$. We have $\Psi(P)=\Psi_{d^*(P)}(P)$, but
now $d^*(P)$ is uniquely defined for all $W$. More generally, $d^*(P)$
represents a uniquely defined optimal rule. Other formulations of the
optimal rule can be obtained by changing the behavior of the rule
$B_0$. Our goal is to construct root-$n$ rate CIs for $\Psi(P_0)$ that
maintain nominal coverage, even at exceptional laws. At non-exceptional
laws, we would like these CIs to belong to and be asymptotically
efficient among the class of regular asymptotially linear (RAL) estimators.

\section{Necessary and sufficient conditions for pathwise
differentiability of \texorpdfstring{$\Psi$}{Psi}} \label{sec:pd}
In this section, we give a necessary and sufficient condition for the
pathwise differentiability of the optimal value parameter $\Psi$. When
it exists, the pathwise derivative in a nonparametric model can be
written as an inner product between an almost surely unique mean zero,
square integrable function known as the canonical gradient and a score
function. The canonical gradient is a key object in nonparametric
statistics. We remind the reader that an estimator $\hat{\Phi}$ is
asymptotically linear for a parameter mapping $\Phi$ at $P_0$ with
influence curve $\mathit{IC}_0$ if
\[
\hat{\Phi}(P_n)-\Phi(P_0)= \frac{1}{n}\sum
_{i=1}^n \mathit{IC}_0(O_i)
+ o_{P_0}\bigl(n^{-1/2}\bigr),
\]
where $E_{0}[\mathit{IC}_0(O)]=0$. The pathwise derivative is important because,
when $\Phi$ is pathwise differentiable in a nonparametric model, any
regular estimator $\hat{\Phi}$ is asymptotically linear with influence
curve $\mathit{IC}_0(O_i)$ if and only if $\mathit{IC}_0$ is the canonical gradient
[\citet
{Bickel1993}]. We discuss negative results for non-pathwise
differentiable parameters and formally define ``regular estimator''
later in this section.

The pathwise derivative of $\Psi$ at $P_0$ can be defined as follows.
Define paths $\{P_\varepsilon:\varepsilon\in\mathbb{R}\}\subset
{\cal M}$
that go through $P_0$ at $\varepsilon=0$, that is, $P_{\varepsilon=0}=P_0$.
In particular, these paths are given by
%
\begin{eqnarray}
\label{eq:fluct} dQ_{W,\varepsilon}&=&\bigl(1+\varepsilon S_W(W)
\bigr)\,dQ_{W,0},
\nonumber
\\
\eqntext{\mbox{where }E_0\bigl[S_W(W)\bigr]=0\mbox{ and }
\displaystyle\sup_w \bigl|S_W(w)\bigr|<\infty;}
\\[-8pt]
\\[-8pt]
\nonumber
dQ_{Y,\varepsilon}(Y|A,W)&=&\bigl(1+\varepsilon S_Y(Y|A,W)
\bigr)\,dQ_{Y,0}(Y|A,W),
\nonumber
\\
\eqntext{\mbox{where }E_0[S_Y|A,W]=0\mbox{ and }\displaystyle\sup
_{w,a,y} \bigl|S_Y(y|a,w)\bigr|<\infty.}
\end{eqnarray}
Above $Q_{W,0}$ and $Q_{Y,0}$ are respectively the marginal
distribution of $W$ and the conditional distribution of $Y$ given $A,W$
under $P_0$. The parameter $\Psi$ is not sensitive to fluctuations of
$g_0(a|w)=\pr(a|w)$, and thus we do not need to fluctuate this portion
of the likelihood. The parameter $\Psi$ is called pathwise
differentiable at $P_0$ if
\[
\frac{d}{d\varepsilon}\Psi(P_\varepsilon)\Big\vert _{\varepsilon=0}=\int
D^*(P_0) (o) \bigl(S_W(w) + S_{Y|A,W}(y|a,w)
\bigr)\,dP_0(o)
\]
for some $P_0$ mean zero, square integrable function $D^*(P_0)$ with
$E_{0}[D^*\times\break (P_0)(O)| A,W]$ almost surely equal to $E_{0}[D^*(P_0)(O)|W]$.
We refer the reader to \citet{Bickel1993} for a more general exposition
of pathwise differentiability.


In \citet{vanderLaanLuedtke14b}, we showed that $\Psi$ is pathwise
differentiable at $P_0$ with canonical gradient $D(d_0^*,P_0)$ if $P_0$
is a non-exceptional law, that is, $\pr(W\notin B_0)=1$. Exceptional
laws were shown to present problems for estimation of optimal rules
indexed by a finite dimensional parameter by \citet{Robins04}, and it
was observed by \citet{RobinsRotnitzky14} that these laws can also
cause problems for unrestricted optimal rules. Here, we show that mean
outcome under the optimal rule is pathwise differentiable under a
slightly more general condition than requiring a non-exceptional law,
namely that
%
\begin{equation}
\pr \Bigl\{w\in\mathcal{W} : w\notin B_0 \mbox{ or } \max
_{a\in
\{0,1\}
} \sigma_0^2(a,w)=0 \Bigr\}=1,
\label{eq:pathcond}
\end{equation}
where $\sigma_0(a,w)\triangleq\sqrt{\Var_{P_0}(Y|A=a,W=w)}$. The
upcoming theorem also gives the converse result, that is, the mean
outcome under the optimal rule is not pathwise differentiable if the
above condition does not hold.

\begin{theorem}\label{thm:notpd}
Assume $\pr(0<g_0(1|W)<1)=1$, $\pr(|Y|<M)=1$ for some $M<\infty$, and
$\Var_{P_0}(D(d_0^*,P_0)(O))<\infty$.
The parameter $\Psi(P_0)$ is pathwise differentiable if and only if
(\ref{eq:pathcond}) holds. If $\Psi$ is pathwise differentiable at
$P_0$, then $\Psi$ has canonical gradient $D(d_0^*,P_0)$ at $P_0$.
\end{theorem}

In the proof of the theorem, we construct fluctuations $S_W$ and $S_Y$
such that
%
\begin{equation}
\lim_{\varepsilon\uparrow0}\frac{\Psi(P_\varepsilon)-\Psi
(P_0)}{\varepsilon}\neq\lim_{\varepsilon\downarrow0}
\frac{\Psi
(P_\varepsilon)-\Psi(P_0)}{\varepsilon} \label{eq:limsnoteq}
\end{equation}
when (\ref{eq:pathcond}) does not hold. It then follows that $\Psi
(P_0)$ is not pathwise differentiable. The left- and right-hand sides
above are referred to as one-sided directional derivatives by \citet
{HiranoPorter2012}.

This condition for the mean outcome differs slightly from that implied
for unrestricted rules in \citet{RobinsRotnitzky14} in that we still
have pathwise differentiability when the $\bar{Q}_{b,0}$ is zero in
some strata but the conditional variance of the outcome given
covariates and treatment is also zero in all of those strata. This
makes sense, given that in this case the blip function could be
estimated perfectly in those strata in any finite sample with treated
and untreated individuals observed in that strata. Though we do not
expect this difference to matter for most data generating distributions
encountered in practice, there are cases where it may be relevant. For
example, if no one in a certain strata is susceptible to a disease
regardless of treatment status, and researchers are unaware of this
{a priori} so that simply excluding this strata from the target
population is not an option, then the treatment effect and conditional
variance are both zero in this strata.

In general, however, we expect that the mean outcome under the optimal
rule will not be pathwise differentiable under exceptional laws
encountered in practice. For this reason, we often refer to
``exceptional laws'' rather than ``laws which do not satisfy (\ref
{eq:pathcond})'' in this work. We do this because the term
``exceptional law'' is well established in the literature, and also
because we believe that there is likely little distinction between
``exceptional law'' and ``laws which do not satisfy (\ref
{eq:pathcond})'' for many problems of interest.

For the definitions of regularity and local unbiasedness, we let
$P_\varepsilon$ be as in (\ref{eq:fluct}), with $g_0$ also fluctuated.
That is, we let $dP_\varepsilon=dQ_{Y,\varepsilon}\times
g_{\varepsilon}\times
dQ_{W,\varepsilon}$, where $g_\varepsilon(A|W)=(1+\varepsilon
S_A(A|W))g_0(A|W)$
with $E_0[S_A(A|W)|W]=0$ and $\sup_{a,w}|S_A(a|w)|<\infty$. The
estimator $\hat{\Phi}$ of $\Phi(P_0)$ is called regular if the
asymptotic distribution of $\sqrt{n}(\hat{\Phi}(P_n)-\Phi(P_0))$ is not
sensitive to small fluctuations in $P_0$. That is, the limiting
distribution of $\sqrt{n}(\hat{\Phi}(P_{n,\varepsilon=1/\sqrt
{n}})-\Phi
(P_{\varepsilon=1/\sqrt{n}}))$ does not depend on $S_W$, $S_A$, or $S_Y$,
where $P_{n,\varepsilon=1/\sqrt{n}}$ is the empirical distribution
$O_1,\ldots,O_n$ drawn i.i.d. from $P_{\varepsilon=1/\sqrt{n}}$. The
estimator $\hat{\Phi}$ is called locally unbiased if the limiting
distribution of $\sqrt{n}(\hat{\Phi}(P_{n,\varepsilon=1/\sqrt
{n}})-\Phi
(P_{\varepsilon=1/\sqrt{n}}))$ has mean zero for all fluctuations $S_W$,
$S_A$ and $S_Y$, and is called asymptotically unbiased (at $P_0$) if
the bias of $\hat{\Phi}(P_n)$ for the parameter $\Phi(P_0)$ is
$o_{P_0}(n^{-1/2})$ at $P_0$.

The non-regularity of a statistical inference problem does not
typically imply the nonexistence of asymptotically unbiased estimators
[see Example~4 of \citet{LiuBrown1993} and the
discussion thereof in \citet{Chen2004}], but rather the non-existence of
\textit{locally} asymptotically unbiased estimators whenever (\ref
{eq:limsnoteq}) holds for some fluctuation [\citet{HiranoPorter2012}].
It is thus not surprising that we are able to find an estimator that is
asymptotically unbiased at a fixed (possibly exceptional) law under
mild assumptions. \citeauthor{HiranoPorter2012} also show that there
does not exist a regular estimator of the optimal value at any law for
which (\ref{eq:limsnoteq}) holds for some fluctuation. That is, no
regular estimators of $\Psi(P_0)$ exist at laws which satisfy the
conditions of Theorem~\ref{thm:notpd} but do not satisfy (\ref
{eq:pathcond}), that is, one must accept the non-regularity of their
estimator when the data is generated from such a law. Note that this
does not rule out the development of locally consistent confidence
bounds similar to those presented by \citet{LaberMurphy11} and \citet
{Laberetal2014}, though such approaches can be conservative when the
estimation problem is non-regular.

In this work, we present an estimator $\hat{\Psi}$ for which $\Gamma_n
\sqrt{n}(\hat{\Psi}(P_n)-\Psi(P_0))$ converges in distribution to a
standard normal distribution for a random standardization term $\Gamma
_n$ under reasonable conditions. Our estimator does not require any
complicated tuning parameters, and thus allows one to easily develop
root-$n$ rate CIs for the optimal value. We show that our estimator is
RAL and efficient at laws which satisfy (\ref{eq:pathcond}) under conditions.\

\section{Inference at exceptional laws} \label{sec:challenge}
\subsection{The challenge}
Before presenting our estimator, we discuss the challenge of estimating
the optimal value at exceptional laws. Suppose $d_n$ is an estimate of
$d_0^*$ and $\hat{\Psi}_{d_n}(P_n)$ is an estimate of $\Psi(P_0)$
relying on the full data set. In \citet{vanderLaanLuedtke2014}, we
presented a targeted minimum loss-based estimator (TMLE) $\hat{\Psi
}_{d_n}(P_n)$ which satisfies
\begin{eqnarray*}
&&\hat{\Psi}_{d_n}(P_n) - \Psi(P_0)\\
&&\qquad=
(P_n-P_0)D\bigl(d_n,P_n^*\bigr)
+ \underbrace {\Psi_{d_n}(P_0)-\Psi(P_0)}_{o_{P_0}(n^{-1/2})\ \mathrm{under\ conditions}}
\mathrel{+} o_{P_0}\bigl(n^{-1/2}\bigr),
\end{eqnarray*}
where we use the notation $Pf = E_P[f(O)]$ for any distribution $P$ and
the second $o_{P_0}(n^{-1/2})$ term is a remainder from a first-order
expansion of $\Psi$. The term $\Psi_{d_n}(P_0)-\Psi(P_0)$ being
$o_{P_0}(n^{-1/2})$ relies on the optimal rule being estimated well in
terms of value and will often prove to be a reasonable condition, even
at exceptional laws (see Theorem~\ref{thm:suffA5} in Section~\ref
{sec:discdnconv}). Here, $P_n^*$ is an estimate of the components of
$P_0$ needed to estimate $D(d_n,P_0)$. To show asymptotic linearity,
one might try to replace $D(d_n,P_n^*)$ with a term that does not rely
on the sample:
\begin{eqnarray*}
(P_n-P_0)D\bigl(d_n,P_n^*
\bigr)&=& (P_n-P_0)D\bigl(d_0^*,P_0
\bigr)
\\
&&{}+ \underbrace{(P_n-P_0) \bigl(D\bigl(d_n,P_n^*
\bigr)-D\bigl(d_0^*,P_0\bigr)\bigr)}_{\mathrm{empirical\ process}}.
\end{eqnarray*}
If $D(d_n,P_n^*)$ belongs to a Donsker class and converges to
$D(d_0^*,P_0)$ in $L^2(P_0)$, then the empirical process term is
$o_{P_0}(n^{-1/2})$ and $\sqrt{n}(\hat{\Psi}_{d_n}(P_n) - \Psi(P_0))$
converges in distribution to a normal random variable with mean zero
and variance $\Var_{P_0}(D(d_0^*,P_0))$ [\citet{vanderVaartWellner1996}].
Note that $D(d_n,P_n^*)$ being consistent for $D(d_0^*,P_0)$ will
typically rely on $d_n$ being consistent for the fixed $d_0^*$ in
$L^2(P_0)$, which we emphasize is \textit{not} implied by $\Psi
_{d_n}(P_0)-\Psi(P_0)=o_{P_0}(n^{-1/2})$. \citet{Zhangetal12b} make this
assumption in the regularity conditions in their Web Appendix A when
they consider an analogous empirical process term in deriving the
standard error of an estimate of the optimal value in a restricted
class. More specifically, \citeauthor{Zhangetal12b} assume a
non-exceptional law and consistent estimation of a fixed optimal rule.
\citet{vanderLaanLuedtke2014} also make such an assumption. If $P_0$
is an exceptional law, then we likely do not expect $d_n$ to be
consistent for any fixed (non-data dependent) function. Rather, we
expect $d_n$ to fluctuate randomly on the set $B_0$, even as the sample
size grows to infinity. In this case, the empirical process term
considered above is not expected to behave as $o_{P_0}(n^{-1/2})$.

Accepting that our estimates of the optimal rule may not stabilize as
sample size grows, we consider an estimation strategy that allows $d_n$
to remain random even as $n\rightarrow\infty$.


\subsection{A thought experiment} \label{sec:flawedarg}
First, we give an erroneous estimation strategy which contains the main
idea of the approach but is not correct in its current form. A
modification is given in the next section. For simplicity, we will
assume that one knows $v_n\triangleq\Var_{P_0}(D(d_n,P_0))$ given an
estimate $d_n$ and, for simplicity, that $v_n$ is almost surely bounded
away from zero. Under reasonable conditions,
\[
v_n^{-1/2} \bigl(\hat{\Psi}_{d_n}(P_n) -
\Psi(P_0) \bigr)= (P_n-P_0)v_n^{-1/2}D
\bigl(d_n,P_n^*\bigr) + o_{P_0}
\bigl(n^{-1/2}\bigr).
\]
The empirical process on the right is difficult to handle because $d_n$
and $v_n$ are random quantities that likely will not stabilize to a
fixed limit at exceptional laws.

As a thought experiment, suppose that we could treat $\{
v_n^{-1/2}D(d_n,P_n^*) : n\}$ as a deterministic sequence, where this
sequence does not necessarily stabilize as sample size grows. In this
case, the Lindeberg--Feller central limit theorem (CLT) for triangular
arrays [see, e.g., \citet{AthreyaLahiri2006}] would allow us to show
that the leading term on the right-hand side converges to a standard
normal random variable. This result relies on inverse weighting by
$\sqrt{v_n}$ so the variance of the terms in the sequence stabilizes to
one as sample size gets large.

Of course, we cannot treat these random quantities as deterministic. In
the next section, we will use the general trick of inverse weighting by
the standard deviation of the terms over which we are taking an
empirical mean, but we will account for the dependence of the estimated
rule $d_n$ on the data by inducing a martingale structure that allows
us to treat a sequence of estimates of the optimal rule as known
(conditional on the past). We can then apply a martingale CLT for
triangular arrays to obtain a limiting distribution for our estimator.

\section{Estimation of and inference for the optimal value} \label
{sec:martcis}
In this section, we present a modified one-step estimator $\hat{\Psi}$
of the optimal value. This estimator relies on estimates of the
treatment mechanism $g_0$, the strata-specific outcome $\bar{Q}_0$, and
the optimal rule $d_0^*$. We first present our estimator, and then
present an asymptotically valid two-sided CI for the optimal value
under conditions. Next, we give conditions under which our estimator is
RAL and efficient, and finally we present a (potentially conservative)
asymptotically valid one-sided CI which lower bounds the mean outcome
under the unknown optimal treatment rule. The one-sided CI uses the
same lower bound from the two-sided CI, but does not require a
condition about the rate at which the value of the optimal rule
converges to the optimal value, or even that the value of the estimated
rule is consistent for the optimal value.

The estimators in this section can be extended to a martingale-based
TMLE for $\Psi(P_0)$. Because the primary purpose of this paper is to
deal with inference at exceptional laws, we will only present an online
one-step estimator and leave the presentation of such a TMLE to future work.

\subsection{Estimator of the optimal value} \label{sec:estimator}
In this section, we present our estimator of the optimal value. Our
procedure first estimates the needed features $g_0$, $\bar{Q}_0$, and
$d_0^*$ of the likelihood based on a small chunk of data, and then
evaluates a one-step estimator with these nuisance function values on
the next chunk of the data. It then estimates the features on the first
two chunks of data, and evaluates the one-step estimator on the next
chunk of data. This procedure iterates until we have a sequence of
estimates of the optimal value. We then output a weighted average of
these chunk-specific estimates as our final estimate of the optimal
value. While the first chunk needs to be large enough to estimate the
desired nuisance parameters, that is, large enough to estimate the
features, all subsequent chunks can be of arbitrary size (as small as a
single observation).

We now formally describe our procedure. Define
\[
\tilde{D}(d,\bar{Q},g) (o)\triangleq\frac{I(a=d(w))}{g(a|w)}\bigl(Y-\bar {Q}(a,w)
\bigr) + \bar{Q}\bigl(d(w),w\bigr).
\]
Let $\{\ell_n\}$ be some sequence of non-negative integers representing
the smallest sample on which the optimal rule is learned. For each
$j=1,\ldots,n$, let $P_{n,j}$ represent the empirical distribution of the
observations $(O_1,O_2,\ldots,O_j)$. Let $g_{n,j}$, $\bar{Q}_{n,j}$, and
$d_{n,j}$ respectively represent estimates of the $g_0$, $\bar{Q}_0$,
and $d_0^*$ based on (some subset of) the observations
$(O_1,\ldots,O_{j-1})$ for all $j>\ell_n$. We subscript each of these
estimates by both $n$ and $j$ because the subsets on which these
estimates are obtained may depend on sample size. We give an example of
a situation where this would be desirable in Section~\ref{sec:nuischunks}.

Define
\[
\tilde{\sigma}_{0,n,j}^2\triangleq\Var_{P_0} \bigl(
\tilde {D}(d_{n,j},\bar{Q}_{n,j},g_{n,j}) (O)
\vert O_1,\ldots,O_{j-1} \bigr).
\]
Let $\tilde{\sigma}_{n,j}^2$ represent an estimate of $\tilde{\sigma
}_{0,n,j}^2$ based on (some subset of) the observations
$(O_1,\ldots,O_{j-1})$. Note that we omit the dependence of $\tilde
{\sigma
}_{n,j}$ and $\tilde{\sigma}_{0,n,j}$ on $d_{n,j}$, $\bar{Q}_{n,j}$,
and $g_{n,j}$ in the notation. Our results apply to any sequence of
estimates $\tilde{\sigma}_{n,j}^2$ which satisfies conditions~(C1) through (C5),
which are stated later in this section. Also define
\[
\Gamma_n\triangleq\frac{1}{n-\ell_n}\sum
_{j=\ell_n+1}^n \tilde {\sigma }_{n,j}^{-1}.
\]
Our estimate $\hat{\Psi}(P_n)$ of $\Psi(P_0)$ is given by
%
\begin{equation}
\label{eq:psihatdef} \hat{\Psi}(P_n)\triangleq\Gamma_n^{-1}
\frac{1}{n-\ell_n}\sum_{j=\ell
_n+1}^n \tilde{
\sigma}_{n,j}^{-1}\tilde{D}_{n,j}(O_j) =
\frac{\sum_{j=\ell_n+1}^n \tilde{\sigma}_{n,j}^{-1}\tilde{D}_{n,j}(O_j)}{\sum_{j=\ell_n+1}^n \tilde{\sigma}_{n,j}^{-1}},
\end{equation}
where $\tilde{D}_{n,j}\triangleq\tilde{D}(d_{n,j},\bar
{Q}_{n,j},g_{n,j})$. We note that the $\Gamma_n^{-1}$ standardization
is used to account for the term-wise inverse weighting so that $\hat
{\Psi}(P_n)$ estimates $\Psi(P_0) = E_0[\tilde{D}(d_0^*,\bar
{Q}_0,g_0)]$. The above looks a lot like a standard augmented inverse
probability weighted (AIPW) estimator, but with $d_0^*$ estimated on
chunks of data increasing in size and with each term in the sum given
weight proportional to an estimate of the conditional variance of that
term. Our estimator constitutes a minor modification of the online
one-step estimator presented in \citet{vanderLaanLendle2014}. In
particular, each term in the sum is inverse weighted by an estimate of
the standard deviation of $\tilde{D}_{n,j}$. For ease of reference, we
will refer to the estimator above as an online one-step estimator.

This estimation scheme differs from sample split estimation, where
features are estimated on half of the data and then a one-step
estimator is evaluated on the remaining half of the data. While one can
show that such estimators achieve valid coverage using Wald-type CIs,
these CIs will generally be approximately $\sqrt{2}$ times larger than
the CIs of our proposed procedure (see the next section) because the
one-step estimator is only applied to half of the data. Alternatively,
one could try averaging two such estimators, where the training and the
one-step sample are swapped between the two estimators. Such a
procedure will fail to yield valid Wald-type CIs due to the
non-regularity of the inference problem: one cannot replace the optimal
rule estimates with their limits because such limits will not generally
exist, and thus the estimator averages over terms with a complicated
dependence structure.

\subsection{Two-sided confidence interval for the optimal value}
Define the remainder terms
\begin{eqnarray*}
R_{1n}&\triangleq&\frac{1}{n-\ell_n}\sum_{j=\ell_n+1}^n
\tilde {\sigma }_{n,j}^{-1}E_{0} \biggl[ \biggl(1-
\frac
{g_0(d_{n,j}(W)|W)}{g_{n,j}(d_{n,j}(W)|W)} \biggr)
\\
&&{}\times \bigl(\bar{Q}_{n,j}\bigl(d_{n,j}(W),W\bigr)-\bar
{Q}_0\bigl(d_{n,j}(W),W\bigr) \bigr) \biggr],
\\
R_{2n}&\triangleq&\frac{1}{n-\ell_n}\sum_{j=\ell_n+1}^n
\frac{\Psi
_{d_{n,j}}(P_0)-\Psi(P_0)}{\tilde{\sigma}_{n,j}}.
\end{eqnarray*}
The upcoming theorem relies on the following assumptions:
\begin{longlist}[(C5)]
\item[(C1)]$n-\ell_n$ diverges to infinity as $n$ diverges to infinity.
%
\item[(C2)] Lindeberg-like condition: for all $\varepsilon>0$, 
\[
\frac{1}{n-\ell_n}\sum_{j=\ell_n+1}^n
E_{0} \biggl[\biggl(\frac
{\tilde{D}_{n,j}(O)}{\tilde{\sigma}_{n,j}} \biggr)^2
T_{n,j}(O)\Big\vert O_1,\ldots,O_{j-1} \biggr]=
o_{P_0}(1),
\]
where $T_{n,j}(O)\triangleq I (\frac{|\tilde
{D}_{n,j}(O)|}{\tilde
{\sigma}_{n,j}}>\varepsilon\sqrt{n-\ell_n} )$.
\item[(C3)]$\frac{1}{n-\ell_n}\sum_{j=\ell_n+1}^{n} \frac{\tilde
{\sigma
}_{0,n,j}^2}{\tilde{\sigma}_{n,j}^2}$ converges to $1$ in probability.
%
\item[(C4)]$R_{1n}=o_{P_0}(n^{-1/2})$. 
\item[(C5)]$R_{2n} = o_{P_0}(n^{-1/2})$. 
\end{longlist}
The assumptions are discussed in Section~\ref{sec:discconds}. We note that
all of our results also hold with $R_{1n}$ and $R_{2n}$ behaving as
$o_{P_0}(1/\sqrt{n-\ell_n})$, though we do not expect this observation
to be of use in practice as we recommend choosing $\ell_n$ so that
$n-\ell_n$ increases at the same rate as $n$.

\begin{theorem} \label{thm:martiptw}
Under {conditions~\textup{(C1)}} through \textup{(C5)}, we have that
\[
\Gamma_n\sqrt{n-\ell_n} \bigl(\hat{\Psi}(P_n)
- \Psi(P_0) \bigr)\rightsquigarrow N(0,1),
\]
where we use ``$\rightsquigarrow$'' to denote convergence in
distribution as the sample size converges to infinity. It follows that
an asymptotically valid $1-\alpha$ CI for $\Psi(P_0)$ is given by
\[
\hat{\Psi}(P_n)\pm z_{1-\alpha/2}\frac{\Gamma_n^{-1}}{\sqrt
{n-\ell_n}},
\]
where $z_{1-\alpha/2}$ denotes the $1-\alpha/2$ quantile of a standard
normal random variable.
\end{theorem}

We have shown that, under very general conditions, the above CI yields
an asymptotically valid $1-\alpha$ CI for $\Psi(P_0)$. We refer the
reader to Section~\ref{sec:discconds} for a detailed discussion of the
conditions of the theorem. We note that our estimator is asymptotically
unbiased, that is, has bias of the order $o_{P_0}(n^{-1/2})$, provided
that $\Gamma_n=O_{P_0}(1)$ and $n-\ell_n$ grows at the same rate as $n$.

Interested readers can consult the proof of {Theorem~\ref
{thm:martiptw}} in the Appendix for a better
understanding of why we proposed the particular estimator given in
{Section~\ref{sec:estimator}}.

\subsection{Conditions for asymptotic efficiency}
We will now show that, if $P_0$ is a non-exceptional law and $d_{n,j}$
has a fixed optimal rule limit $d_0$, then our online estimator is RAL
for $\Psi(P_0)$. The upcoming corollary makes use of the following
consistency conditions for some fixed rule $d_0$ which falls in the
class of optimal rules given in (\ref{eq:optclass}):
%
\begin{eqnarray}
&&\frac{1}{n-\ell_n}\sum_{j=\ell_n+1}^n
E_{0} \bigl[\bigl(d_{n,j}(W)-d_0(W)
\bigr)^2\vert O_1,\ldots,O_{j-1} \bigr]=
o_{P_0}(1), \label{eq:dconscor}
\\
\label{eq:Qconscor}&&\frac{1}{n-\ell_n}\sum_{j=\ell_n+1}^n
E_{0} \bigl[ \bigl(\bar {Q}_{n,j}\bigl(d_0(W),W
\bigr)-\bar{Q}_0\bigl(d_0(W),W\bigr)
\bigr)^2\vert O_1,\ldots,O_{j-1} \bigr]
\nonumber
\\[-8pt]
\\[-8pt]
\nonumber
&&\qquad=
o_{P_0}(1),
\\
\label{eq:gconscor}&&\frac{1}{n-\ell_n}\sum_{j=\ell_n+1}^n
E_{0} \bigl[ \bigl(g_{n,j}\bigl(d_0(W)|W
\bigr)-g_0\bigl(d_0(W)|W\bigr) \bigr)^2
\vert O_1,\ldots ,O_{j-1} \bigr]
\nonumber
\\[-8pt]
\\[-8pt]
\nonumber
&&\qquad= o_{P_0}(1).
\end{eqnarray}
It also makes use of the following conditions, which are, respectively,
slightly stronger than {conditions~(C1)} and (C3):
\begin{longlist}[(C1$'$)]
\item[(C1$'$)]$\ell_n=o(n)$. 
\item[(C3$'$)]$\frac{1}{n-\ell_n}\sum_{j=\ell_n+1}^{n} \llvert \frac{\tilde
{\sigma
}_{0,n,j}^2}{\tilde{\sigma}_{n,j}^2}-1\rrvert \rightarrow0$ in
probability. 
\end{longlist}

\begin{corollary} \label{cor:ral}
Suppose that {conditions~\textup{(C1$'$)}}, \textup{(C2)},
\textup{(C3$'$)},
\textup{(C4)} and \textup{(C5)} hold. Also suppose that $\pr
(\delta<g_0(1|W)<1-\delta)=1$ for some $\delta>0$, the estimates
$g_{n,j}$ are bounded away from zero with probability $1$, $Y$ is
bounded, the estimates $\bar{Q}_{n,j}$ are uniformly bounded,\vspace*{1pt} $\ell
_n=o(n)$, and that, for some fixed optimal rule $d_0$, (\ref
{eq:dconscor}), (\ref{eq:Qconscor}) and (\ref{eq:gconscor}) hold.
Finally, assume that $\Var_{P_0}(\tilde{D}(d_0,\bar{Q}_0,g_0))>0$ and
that, for some $\delta_0>0$, we have that
\[
\pr \Bigl(\inf_{j,n}\tilde{\sigma}_{n,j}^2>
\delta_0 \Bigr)=1,
\]
where the infimum is over natural number pairs $(j,n)$ for which $\ell
_n<j\le n$. Then we have that
%
\begin{equation}
\Gamma_n^{-1}\rightarrow\Var_{P_0}\bigl(
\tilde{D}(d_0,\bar {Q}_0,g_0)\bigr)\qquad\mbox{in probability as }n\rightarrow\infty. \label{eq:sameci}
\end{equation}
Additionally,
%
\begin{equation}
\hat{\Psi}(P_n) - \Psi(P_0)= \frac{1}{n}\sum
_{i=1}^n D(d_0,P_0)
+ o_{P_0}(1/\sqrt{n}). \label{eq:raleff}
\end{equation}
That is, $\hat{\Psi}(P_n)$ is asymptotically linear with influence
curve $D(d_0,P_0)$. Under the conditions of this corollary, it follows
that $P_0$ satisfies (\ref{eq:pathcond}) if and only if $\hat{\Psi
}(P_n)$ is RAL and asymptotically efficient among all such RAL estimators.
\end{corollary}

We note that (\ref{eq:sameci}) combined with (C1$'$)
implies that the CI given in Theorem~\ref{thm:martiptw} asymptotically has
the same width [up to an $o_{P_0}(n^{-1/2})$ term] as the CI which
treats (\ref{eq:raleff}) and $D(d_0,P_0)$ as known and establishes a
typical Wald-type CI about $\hat{\Psi}(P_n)$.

The empirical averages over $j$ in (\ref{eq:dconscor}), (\ref
{eq:Qconscor}) and (\ref{eq:gconscor}) can easily be dealt with using
{Lemma~\ref{lem:rnoprootn}}, presented in
Section~\ref{sec:discsigmacons}. Essentially, we have required that
$d_{n,j}$, $\bar{Q}_{n,j}$ and $g_{n,j}$ are consistent for $d_0$,
$\bar
{Q}_0$ and $g_0$ as $n$ and $j$ get large, where $d_0$ is some fixed
optimal rule. One would expect such a fixed limiting rule $d_0$ to
exist at a non-exceptional law for which the optimal rule is (almost
surely) unique. If $g_0$ is known, then we do not need $\bar{Q}_{n,j}$
to be consistent for $\bar{Q}_0$ to get asymptotic linearity, but
rather that $\bar{Q}_{n,j}$ converges to some possibly misspecified
fixed limit~$\bar{Q}$. 

\subsection{Lower bound for the optimal value}
It would likely be useful to have a conservative lower bound on the
optimal value in practice. If policymakers were to implement an optimal
individualized treatment rule whenever the overall benefit is greater
than some fixed threshold, that is, $\Psi(P_0)>v$ for some fixed $v$,
then a one-sided CI for $\Psi(P_0)$ would help facilitate the decision
to implement an individualized treatment strategy in the population.

The upcoming theorem shows that the lower bound from the $1-2\alpha$ CI
yields a (potentially conservative) asymptotic $1-\alpha$ CI for the
optimal value. If $d_0^*$ is estimated well in the sense of
{condition~(C5)}, then the asymptotic coverage
is exact. Define
\[
{LB}_n(\alpha)\triangleq\hat{\Psi}(P_n)-
z_{1-\alpha}\frac
{\Gamma
_n^{-1}}{\sqrt{n-\ell_n}}.
\]

\begin{theorem} \label{thm:cilb}
Under {conditions~\textup{(C1)}} through \textup{(C4)}, we have that
\[
\liminf_{n\rightarrow\infty} \pr \bigl(\Psi(P_0)>{LB}_n(
\alpha ) \bigr)\ge1-\alpha.
\]
If {condition~\textup{(C5)}} also holds, then
\[
\lim_{n\rightarrow\infty} \pr \bigl(\Psi(P_0)>{LB}_n(
\alpha ) \bigr)= 1-\alpha.
\]
\end{theorem}

The above condition should not be surprising, as we base our CI for
$\Psi(P_0)$ on a weighted combination of estimates of $\Psi
_{d_{n,j}}(P_0)$ for $j<n$. Because $\Psi(P_0)\ge\Psi_{d_{n,j}}(P_0)$
for all such $j$, we would expect that the lower bound of the $1-\alpha
$ CI given in the previous section provides a valid $1-\alpha/2$
one-sided CI for $\Psi(P_0)$. Indeed this is precisely what we see in
the proof of the above theorem.

\subsection{Coverage for the value of the rule estimated on the entire
data set}
Suppose one wishes to evaluate the coverage of our CI for the data
dependent parameter $\Psi_{d_n}(P_0)$, where $d_n$ is an estimate of
the optimal rule based on the entire data set of size $n$. We make two
key assumptions in this section, namely that there exists some real
number $\psi_1$ such that:
\begin{longlist}[(C6)]
\item[(C6)]$\Gamma_n(\Psi_{d_n}(P_0)-\psi_1)= o_{P_0}(n^{-1/2})$. 
\item[(C7)]$\frac{1}{n-\ell_n}\sum_{j=\ell_n+1}^n \frac{\Psi
_{d_{n,j}}(P_0)-\psi_1}{\tilde{\sigma}_{n,j}}= o_{P_0}(n^{-1/2})$.
\end{longlist}
Typically, $\Gamma_n=O_{P_0}(1)$ so that {condition~(C6)} is the same as $\Psi
_{d_n}(P_0)=\psi_1+o_{P_0}(1/\sqrt{n})$. As will become apparent after
reading {Section~\ref{sec:discconds}},
{condition~(C6)} will typically
imply (C7) (see {Lemma~\ref{lem:rnoprootn}}).
{Theorem~\ref{thm:suffA5}}
shows that {condition (C6)} is
often reasonable with $\psi_1=\Psi(P_0)$, though we do not require that
$\psi_1=\Psi(P_0)$.

\begin{theorem} \label{thm:daparam}
Suppose {conditions~\textup{(C1)}} through
\textup{(C4)} and {conditions~\textup{(C6)}} and \textup{(C7)} hold. Then
\[
\Gamma_n \sqrt{n-\ell_n} \bigl(\hat{
\Psi}(P_n)-\Psi _{d_n}(P_0) \bigr)
\rightsquigarrow N(0,1).
\]
Thus, the same CI given in {Theorem~\ref{thm:martiptw}} is an
asymptotically valid $1-\alpha$ CI for $\Psi_{d_n}(P_0)$.
\end{theorem}


\section{Computationally efficient estimation schemes} \label{sec:compeff}
Computing $\hat{\Psi}(P_n)$ may initially seem computationally
demanding. In this section, we discuss two estimation schemes which
yield computationally simple routines. 

\subsection{Computing the features on large chunks of the data} \label
{sec:nuischunks}
One can compute the estimates of $\bar{Q}_0$, $g_0$ and $d_0^*$ far
fewer than $n-\ell_n$ times. For each $j$, the estimates $\bar
{Q}_{n,j}$, $g_{n,j}$, and $d_{n,j}$ may rely on any subset of the
observations $O_1,\ldots,O_{j-1}$. Thus, one can compute these estimators
on $S$ increasing subsets of the data, where the first subset consists
of observations $O_1,\ldots,O_{\ell_n}$ and each of the $S-1$ remaining
samples adds a $1/S$ proportion of the remaining $n-\ell_n$
observations. Note that this scheme makes use of the fact that, for
fixed $j$, the feature estimates, indexed by $n$ and $j$, for example,
$d_{n,j}$, may rely on different subsets of observations
$O_1,\ldots,O_{j-1}$ for different sample sizes $n$. 
\subsection{Online learning of the optimal value}
Our estimator was inspired by online estimators which can operate on
large data sets that will not fit into memory. These estimators use
online prediction and regression algorithms which update the initial
fit based on previously observed estimates using new observations as
they are read into memory. Online estimators of pathwise differentiable
parameters were introduced in \citet{vanderLaanLendle2014}. Such
estimation procedures often require estimates of features of the
likelihood, which can be obtained using modern online regression and
classification approaches [see,
e.g., \citet{Zhang2004,Langfordetal2009,Lutsetal2013}]. Our estimator
constitutes a slight modification of the one-step online estimator
presented by \citet{vanderLaanLendle2014}, and thus all discussion of
computational efficiency given in that paper applies to our case.

For our estimator, one could use online estimators of $\bar{Q}_0$,
$g_0$ and $d_0^*$, and then update these estimators as the index $j$ in
the sum in (\ref{eq:psihatdef}) increases. Calculating the standard
error estimate $\tilde{\sigma}_{n,j}$ will typically require access to
an increasing subset of the past observations, that is, as sample size
grows one may need to hold a growing number of observations in memory.
If one uses a sample standard deviation to estimate $\tilde{\sigma
}_{0,n,j}$ based on subset of observations $O_1,\ldots,O_{j-1}$, the
results we present in Section~\ref{sec:discsigmacons} will indicate that
one really only needs that the number of points on which $\tilde
{\sigma
}_{0,n,j}$ is estimated grows with $j$ rather than at the same rate as
$j$. This suggest that, if computation time or system memory is a
concern for calculating $\tilde{\sigma}_{n,j}$, then one could
calculate $\tilde{\sigma}_{n,j}$ based on some $o(j)$ subset of
observations $O_1,\ldots,O_{j-1}$. 


\section{Discussion of the conditions of Theorem~\texorpdfstring{\protect\ref
{thm:martiptw}}{2}} \label
{sec:discconds}
For ease of notation, we will assume that, for all $j>\ell_n$, we do
not modify our feature estimates based on the first $j-1$ data points
as the sample size grows. That is, for all sample sizes $m,n$ and all
$j\le\min\{m,n\}$, $d_{n,j} = d_{m,j}$, $\bar{Q}_{n,j} = \bar
{Q}_{m,j}$, $g_{n,j}=g_{m,j}$, and $\tilde{\sigma}_{n,j}=\tilde
{\sigma
}_{m,j}$. One can easily extend all of the discussion in this section
to a more general case where, for example, $d_{n,j}\neq d_{m,j}$ for
$n\neq m$. This may be useful if the optimal rule is estimated in
chunks of increasing size as was discussed in Section~\ref{sec:nuischunks}.
To make these object's lack of dependence on $n$ clear, in this section
we will denote $d_{n,j}$, $\bar{Q}_{n,j}$, $g_{n,j}$, $\tilde{\sigma
}_{n,j}$, and $\tilde{\sigma}_{0,n,j}$ as $d_{j}$, $\bar{Q}_{j}$,
$g_j$, $\tilde{\sigma}_j$ and $\tilde{\sigma}_{0,j}$. This will also
help make it clear when $o_{P_0}$ notation refers to behavior as $j$,
rather than $n$, goes to infinity.

For our discussion, we assume there exists a (possibly unknown) $\delta
_0>0$ such that
%
\begin{equation}
\pr \Bigl(\inf_{j>\ell_n}\tilde{\sigma}_{0,j}^2>
\delta_0 \Bigr)=1, \label
{eq:condvarlb}
\end{equation}
where the probability statement is over the i.i.d. draws
$O_1,O_2,\ldots.$
The above condition is not necessary, but will make our discussion of
the conditions more straightforward. 

\subsection{Discussion of {condition~\texorpdfstring{\textup{(C1)}}{(C1)}}}
\label{eq:discelln}
We cannot apply the martingale CLT in the proof of Theorem~\ref
{thm:martiptw} if $n-\ell_n$ does not grow with sample size.
Essentially, this condition requires that a non-negligible proportion
of the data is used to actually estimate the mean outcome under the
optimal rule. One option is to have $n-\ell_n$ grow at the same rate as
$n$ grows, which holds if, for example, $\ell_n=pn$ for some fixed
proportion $p$ of the data. This allows our CIs to shrink at a root-$n$
rate. One might prefer to have $\ell_n=o(n)$ so that $\frac{n-\ell
_n}{n}$ converges to $1$ as sample size grows. In this case, we can
show that our estimator is asymptotically linear and efficient at
non-exceptional laws under conditions, as we did in {Corollary~\ref{cor:ral}}.

\subsection{Discussion of {condition~\texorpdfstring{\textup{(C2)}}{(C2)}}}
\label{eq:discLind}
This is a standard condition that yields a martingale CLT for
triangular arrays [\citet{Gaenssleretal1978}]. The condition ensures that
the variables which are being averaged have sufficiently thin tails.
While it is worth stating the condition in general, it is easy to
verify that the condition is implied by the following three more
straightforward conditions:
\begin{itemize}
\item(\ref{eq:condvarlb}) holds.
\item$Y$ is a bounded random variable.
\item There exists some $\delta>0$ such that $\pr(\delta< g_j(1|W)
<1-\delta)=1$ with probability $1$ for all $j$.
\end{itemize}
Indeed, under the latter two conditions $|\tilde{D}_{n,j}(O)|<C$ is
almost surely bounded for some $C>0$, and thus (\ref{eq:condvarlb})
yields that $|\tilde{D}_{n,j}(O)\tilde{\sigma}_{n,j}^{-1}|<C\delta
_0^{-1}<\infty$ with probability $1$. For all $\varepsilon>0$,
$\varepsilon
\sqrt{n-\ell_n}>C\delta_0^{-1}$ for all $n$ large enough under
condition~(C1). Thus, $T_{n,j}$ from condition~(C2) is equal to zero with probability $1$ for all $n$ large enough.

\subsection{Discussion of {condition~\texorpdfstring{\textup{(C3)}}{(C3)}}} \label{sec:discsigmacons}
This is a rather weak condition given that $\tilde{\sigma}_{0,j}$ still
treats $d_{j}$ as random. Thus, this condition does not require that
$d_{j}$ stabilizes as $j$ gets large. Suppose that
%
\begin{equation}
\tilde{\sigma}_j^2-\tilde{\sigma}_{0,j}^2
= o_{P_0}(1). \label{eq:sigcons}
\end{equation}
By (\ref{eq:condvarlb}) and the continuous mapping theorem, it follows that
%
\begin{equation}
\frac{\tilde{\sigma}_{0,j}^{2}}{\tilde{\sigma}_j^2}-1=o_{P_0}(1). \label
{eq:siginvcons}
\end{equation}
The following general lemma will be useful in establishing
{conditions~(C3)}, (C4) and (C5).

\begin{lemma} \label{lem:rnoprootn}
Suppose that $R_j$ is some sequence of (finite) real-valued random
variables such that $R_j=o_{P_0}(j^{-\beta})$ for some $\beta\in[0,1)$,
where we assume that each $R_j$ is measurable with respect to the
sigma-algebra generated by $(O_1,\ldots,O_j)$. Then
\[
\frac{1}{n}\sum_{j=1}^n
R_j = o_{P_0}\bigl(n^{-\beta}\bigr).
\]
\end{lemma}

Applying the above lemma with $\beta=0$ to (\ref{eq:siginvcons}) shows
that {condition~\textup{(C3)}}
holds provided that (\ref{eq:condvarlb}) and (\ref{eq:sigcons}) hold.
We will use the above lemma with $\beta=1/2$ when discussing
{conditions~\textup{(C4)}} and \textup{(C5)}.

It remains to show that we can construct a sequence of estimators such
that (\ref{eq:sigcons}) holds. Suppose we estimate $\tilde{\sigma
}_{0,j}^2$ with
%
\begin{equation}
\tilde{\sigma}_{j}^2\triangleq\max \Biggl\{
\delta_j, \frac
{1}{j-1}\sum_{i=1}^{j-1}
\tilde{D}_{j}^2(O_i)- \Biggl(
\frac{1}{j-1}\sum_{i=1}^{j-1}
\tilde{D}_{j}(O_i) \Biggr)^2 \Biggr\},
\label{eq:sig0bdd}
\end{equation}
where $\{\delta_j\}$ is a sequence that may rely on $j$ and each
$\tilde
{D}_{n,j}=\tilde{D}_j$ for all $n\ge j$. We use $\delta_j$ to ensure
that $\tilde{\sigma}_{j}^{-2}$ is well defined (and finite) for all
$j$. If a lower bound $\delta_0$ on $\tilde{\sigma}_{0,j}^2$ is known
then one can take $\delta_j=\delta_0$ for all $j$. Otherwise, one can
let $\{\delta_j\}$ be some sequence such that $\delta_j\downarrow0$ as
$j\rightarrow\infty$.

Note that $\tilde{\sigma}_{j}^2$ is an empirical process because it
involves sums over observations $O_1,\ldots,O_{j-1}$, and functions
$\tilde
{D}_{j}$ which were estimated on those same observations. The following
theorem gives sufficient conditions for (\ref{eq:sigcons}), and thus
{condition~(C3)}, to hold.

\begin{theorem} \label{thm:iptwsigmaconssat}
Suppose (\ref{eq:condvarlb}) holds and that $ \{\tilde{D}(d,\bar
{Q},g) : d,\bar{Q},g \}$ is a $P_0$ Glivenko--Cantelli (GC) class
with an integrable envelope function, where $d$, $\bar{Q}$ and $g$ are
allowed to vary over the range of the estimators of $d_0^*$, $\bar
{Q}_0$, and $g_0$. Let $\tilde{\sigma}_j^2$ be defined as in (\ref
{eq:sig0bdd}). Then we have that $\tilde{\sigma}_j^2-\tilde{\sigma
}_{0,j}^2= o_{P_0}(1)$. It follows that (\ref{eq:siginvcons}) and
{condition~\textup{(C3)}} are satisfied.
\end{theorem}

We thus only make the very mild assumption that our estimators of
$d_0^*$, $\bar{Q}_0$ and $g_0$ belong to GC classes. Note that this
assumption is much milder than the typical Donsker condition needed
when attempting to establish the asymptotic normality of a (non-online)
one-step estimator. An easy sufficient condition for a class to have a
finite envelope function is that it is uniformly bounded, which occurs
if the conditions discussed in Section~\ref{eq:discLind} hold.

\subsection{Discussion of {condition~\texorpdfstring{\textup{(C4)}}{(C4)}}}
This condition is a weighted version of the typical double robust
remainder appearing in the analysis of the AIPW estimator. Suppose that
%
\begin{eqnarray}\label{eq:R1prime}
&&E_{0} \biggl[ \biggl(1-\frac{g_0(d_{j}(W)|W)}{g_{j}(d_{j}(W)|W)} \biggr) \bigl(
\bar{Q}_{j}\bigl(d_{j}(W),W\bigr)-\bar{Q}_0
\bigl(d_{j}(W),W\bigr) \bigr) \biggr]
\nonumber
\\[-8pt]
\\[-8pt]
\nonumber
&&\qquad=o_{P_0}
\bigl(j^{-1/2}\bigr).
\end{eqnarray}
If $g_0$ is known (as in an RCT without missingness) and one takes each
$g_j=g_0$, then the above ceases to be a condition as the left-hand
side is always zero. We note that the only condition on $\bar{Q}_j$
appears in {condition~(C4)}, so that
if $R_{1n}=0$ as in an RCT without missingness then we do not require
that $\bar{Q}_j$ stabilizes as $j$ grows. A typical AIPW estimator
require the estimate of $\bar{Q}_0$ to stabilize as sample size grows
to get valid inference, but here we have avoided this condition in the
case where $g_0$ is known by using the martingale structure and inverse
weighting by the standard error of each term in the definition of $\hat
{\Psi}(P_n)$.

More generally, {Lemma~\ref{lem:rnoprootn}}
shows that {condition~(C4)} holds if
(\ref{eq:siginvcons}) and (\ref{eq:R1prime}) hold and $\pr
(0<g_j(1|W)<1)=1$ with probability $1$ for all $j$. One can apply the
Cauchy--Schwarz inequality and take the maximum over treatment
assignments to see that (\ref{eq:R1prime}) holds if
\[
\max \biggl\{\frac{\Vert g_j(a|W)-g_0(a|W)\Vert_{2,P_0} \Vert\bar
{Q}_j(a,W)-\bar{Q}_0(a,W)\Vert_{2,P_0}}{g_j(a|W)} : a=0,1 \biggr\}
\]
is $o_{P_0}(j^{-1/2})$. If $g_0$ is not known, the above shows that
then (\ref{eq:R1prime}) holds if $g_0$ and $\bar{Q}_0$ are estimated well.

\subsection{Discussion of {condition~\texorpdfstring{\textup{(C5)}}{(C5)}}}
\label{sec:discdnconv}
This condition requires that we can estimate $d_0^*$ well as sample
size gets large. We now give a theorem which will help us to establish
{condition~(C5)} under reasonable
conditions. The theorem assumes the following margin assumption: for
some $\alpha>0$,
%
\begin{equation}
\pr\bigl(0<\bigl|\bar{Q}_{b,0}(W)\bigr|\le t\bigr)\lesssim t^{\alpha}\qquad
\forall t>0, \label{eq:ma}
\end{equation}
where ``$\lesssim$'' denotes less than or equal to up to a nonnegative
constant. This assumption is a direct restatement of Assumption (MA)
from \citet{AudibertTsybakov2007} and was considered earlier by \citet
{Tsybakov04}. Note that this theorem is similar in spirit to Lemma~1 in
\citet{vanderLaanLuedtke2014}, but relies on weaker, and we believe
more interpretable, assumptions.

\begin{theorem} \label{thm:suffA5}
Suppose (\ref{eq:ma}) holds for some $\alpha>0$ and that we have an
estimate $\bar{Q}_{b,n}$ of $\bar{Q}_{b,0}$ based on a sample of size
$n$. If $\Vert\bar{Q}_{b,n}-\bar{Q}_{b,0}\Vert_{2,P_0}=o_{P_0}(1)$, then
\[
\bigl|\Psi_{d_n}(P_0)-\Psi_{d_0^*}(P_0)\bigr|
\lesssim\Vert\bar {Q}_{b,n}-\bar {Q}_{b,0}\Vert_{2,P_0}^{2(1+\alpha)/(2+\alpha)},
\]
where $d_n$ is the function $w\mapsto I(\bar{Q}_{b,n}(w)>0)$. If
$\Vert\bar{Q}_{b,n}-\bar{Q}_{b,0}\Vert_{\infty,P_0}=o_{P_0}(1)$, then
\begin{eqnarray*}
&&\bigl|\Psi_{d_n}(P_0)-\Psi_{d_0^*}(P_0)\bigr|
\\
&&\qquad\le\Vert\bar{Q}_{b,n}-\bar{Q}_{b,0}\Vert_{\infty,P_0}\pr
\bigl(0<\bigl\llvert \bar {Q}_{b,0}(W)\bigr\rrvert \le\Vert
\bar{Q}_{b,n}-\bar{Q}_{b,0}\Vert _{\infty,P_0} \bigr)
\\
&&\qquad\lesssim\Vert\bar{Q}_{b,n}-\bar{Q}_{b,0}\Vert_{\infty
,P_0}^{1+\alpha}.
\end{eqnarray*}
\end{theorem}

The above theorem thus shows that $\Psi_{d_j}(P_0)-\Psi
_{d_0^*}(P_0)=o_{P_0}(j^{-1/2})$ the distribution of $|\bar
{Q}_{b,0}(W)|$ and our estimates of $\bar{Q}_{b,0}$ satisfy reasonable
conditions. If additionally $\tilde{\sigma}_{0,j}$ is estimated well in
the sense of (\ref{eq:siginvcons}), then an application of {Lemma~\ref
{lem:rnoprootn}} shows that {condition~(C5)} is satisfied.

The first part of the proof of Theorem~\ref{thm:suffA5} is essentially a
restatement of Lemma~5.2 in \citet{AudibertTsybakov2007}. Figure A.1 in
Supplementary Appendix C shows various densities which satisfy (\ref
{eq:ma}) at different values of $\alpha$, and also the slowest rate of
convergence for the blip function estimates for which Theorem~\ref
{thm:suffA5} implies {condition~(C5)}. As illustrated in
the figure, $\alpha>1$ implies that
$p_{b,0}(t)\rightarrow0$ as $t\rightarrow0$. Given that we are
interested in laws where $\pr(\bar{Q}_{b,0}(W)=0)>0$, it is unclear how
likely we are to have that $\alpha>1$ when $W$ contains only continuous
covariates. One might, however, believe that the density is bounded
near zero so that (\ref{eq:ma}) is satisfied at $\alpha=1$.


If $\Vert\bar{Q}_{b,n}-\bar{Q}_{b,0}\Vert_{\infty
,P_0}=o_{P_0}(1)$, then the above
theorem indicates an arbitrarily\vadjust{\goodbreak} fast rate for $\Psi_{d_n}(P_0)-\Psi
_{d_0^*}(P_0)$ when there is a margin around zero, that is, $\pr
(0<|\bar
{Q}_{b,0}(W)|\le t)=0$ for some $t>0$. In fact, $\Psi_{d_n}(P_0)-\Psi
_{d_0^*}(P_0)=0$ with probability approaching $1$ in this case. Such a
margin will exist when $W$ is discrete.

One does not have to use a plug-in estimator for the blip function to
estimate the mean outcome under the optimal rule. One could also use
one of the weighted classification approaches, often known as outcome
weighted learning (OWL), recently discussed in the literature to
estimate the optimal rule [\citet
{QianMurphy11,Zhaoetal12,Zhangetal12,RubinvanderLaan12}]. In some
cases, we expect these approaches to give better estimates of the
optimal rule than methods which estimate the conditional outcomes, so
using them may make {condition~(C5)}
more plausible. In \citet{LuedtkevanderLaan2014}, we describe an
ensemble learner that can combine estimators from both the Q-learning
and weighted classification frameworks.

\section{Simulation methods} \label{sec:sims}
We ran four simulations. Simulation D-E is a point treatment case,
where the treatment may rely on a single categorical covariate $W$.
Simulations C-NE and C-E are two different point treatment simulations
where the treatment may rely on a single continuous covariate $W$.
Simulation C-NE uses a non-exceptional law, while simulation C-E uses
an exceptional law. Simulation TTP-E gives simulation results for a
modification of the two time point treatment simulation presented by
\citet{vanderLaanLuedtke2014}, where the data generating distribution
has been modified so the second time point treatment has no effect on
the outcome. This simulation uses the extension to multiple time point
treatments given in Supplementary Appendix B [\citet{LuedtkevanderLaan2014b}].

Each simulation setting was run over 2000 Monte Carlo draws to evaluate
the performance of our new martingale-based method and a classical (and
for exceptional laws incorrect) one-step estimator with Wald-type CIs.
{Table~\ref{tab:nelln}} shows the combinations of
sample size ($n$) and initial chunk size ($\ell_n$) considered for each
estimator. All simulations were run in \texttt{R} [\citet{R2014}].\vspace*{6pt}

\begin{table}[b]
\caption{Primary combinations of sample size ($n$) and initial chunk
size ($\ell_n$) considered in each simulation. Different choices of
$\ell_n$ were considered for C-NE and C-E to explore the sensitivity of
the estimator to the choice of $\ell_n$}
\label{tab:nelln}
\begin{tabular*}{\textwidth}{@{\extracolsep{\fill}}lc@{}}
\hline
\textbf{Simulation} & $\bolds{(n,\ell_n)}$ \\
\hline
D-E & $(1000,100)$, $(4000,100)$ \\
C-NE, C-E, TTP-E & $(250,25)$, $(1000,25)$, $(4000,100)$ \\
\hline
\end{tabular*}
\end{table}

\subsection{Simulation D-E: Discrete $W$} \label{sec:catsim}
\subsection*{Data}
This simulation uses a discrete baseline covariate $W$ with four
levels, a dichotomous treatment $A$, and a binary outcome $Y$. The data
is generated by drawing i.i.d. samples as follows:
\begin{eqnarray*}
W&\sim&\operatorname{Uniform}\{0,1,2,3\},
\\
A|W&\sim&\operatorname{Binomial}(0.5 + 0.1W),
\\
Y|A,W&\sim&\operatorname{Binomial} \bigl(0.4 + 0.2I(A=1,W=0) \bigr),
\end{eqnarray*}
where $\operatorname{Uniform}\{0,1,2,3\}$ is the discrete distribution which
returns each of $0$, $1$, $2$ and $3$ with probability $1/4$. The above
is an exceptional law because $\bar{Q}_{b,0}(w)=0$ for $w\neq0$. The
optimal value is $0.45$.

\subsection*{Estimation methods}
For each $j=\ell_n+1,\ldots,n$, we used the nonparametric maximum
likelihood estimator generated by the first $j-1$ samples to estimate
$P_0$ and the corresponding plug-in estimators to estimate all of the
needed features of the likelihood, including the optimal rule. We used
the sample standard deviation of $\tilde{D}_{n,j}(O_1),\ldots,\tilde
{D}_{n,j}(O_{j-1})$ to estimate $\tilde{\sigma}_{0,j}$.

\subsection{Simulations C-NE and C-E: Continuous univariate $W$}
\label
{sec:univsim}
\subsubsection*{Data}
This simulation uses a single continuous baseline covariate $W$ and
dichotomous treatment $A$ which are sampled as follows:
\begin{eqnarray*}
W&\sim&\operatorname{Uniform}(-1,1),
\\
A|W&\sim&\operatorname{Binomial}(0.5 + 0.1W).
\end{eqnarray*}
We consider two distributions for the binary outcome $Y$. The first
distribution (C-NE) is a non-exceptional law with $Y|A,W$ drawn from to
a $\operatorname{Binomial}(\bar{Q}_0^{\mathrm{n\mbox{-}e}}(A, W))$, where
\begin{eqnarray*}
\bar{Q}_0^{\mathrm{n\mbox{-}e}}(A,W)-\frac{3}{10}\triangleq
\cases{ -W^3 + W^2 -\frac{1}{3}W +
\frac{1}{27},&\quad $\mbox{if }A=1\mbox{ and }W\ge0$, \vspace*{2pt}
\cr
\frac{3}{4}W^3 + W^2 -\frac{1}{3}W +
\frac{1}{27},&\quad $\mbox{if }A=1\mbox { and }W< 0,$ \vspace*{2pt}
\cr
0,&\quad $\mbox{if }A=0$.} %
\end{eqnarray*}
The optimal value of approximately $0.388$ was estimated using $10^8$
Monte Carlo draws. The second distribution (C-E) is an exceptional law
with $Y|A,W$ drawn from to a $\operatorname{Binomial}(\bar{Q}_0^{\mathrm
{e}}(A,W))$, where for $\tilde{W}\triangleq W+5/6$ we define
\begin{eqnarray*}
\bar{Q}_0^{\mathrm{e}}(A,W)-\frac{3}{10}\triangleq
\cases{-\tilde{W}^3 + \tilde{W}^2 -
\frac{1}{3}\tilde{W} + \frac
{1}{27},&\quad $\mbox{if }A=1\mbox{ and }W<
-1/2,$ \vspace*{2pt}
\cr
-W^3 + W^2 -\frac{1}{3}W +
\frac{1}{27},&\quad $\mbox{if }A=1\mbox{ and }W> 1/3,$ \vspace*{2pt}
\cr
0,&\quad $\mbox{otherwise.}$} %
\end{eqnarray*}
The above distribution is an exceptional law because $\bar
{Q}_0^{\mathrm
{e}}(1,w)-\bar{Q}_0^{\mathrm{e}}(0,w)=0$ whenever $w\in [-\frac
{1}{2},\frac{1}{3} ]$. The optimal value of approximately $0.308$
was estimated using $10^8$ Monte Carlo draws.

\subsubsection*{Estimation methods}

To show the flexibility of our estimation procedure with respect to
estimators of the optimal rule, we estimated the blip functions using a
Nadaraya--Watson estimator, where we behave as though $g_0$ is unknown
when computing the kernel estimate. For the next simulation setting, we
use the ensemble learner from \citet{LuedtkevanderLaan2014} that we
suggest using in practice. Here, we estimated
\[
\bar{Q}_{b,n}^h (w)\triangleq\frac{\sum_{i=1}^n y_i a_i K
({(w-w_i)}/{h} )}{\sum_{i=1}^n a_i K ({(w-w_i)}/{h} )} -
\frac{\sum_{i=1}^n y_i (1-a_i)K ({(w-w_i)}/{h} )}{\sum_{i=1}^n (1-a_i)K ({(w-w_i)}/{h} )},
\]
where $K(u)\triangleq\frac{3}{4}(1-u^2)I(|u|\le1)$ is the
Epanechnikov kernel and $h$ is the bandwidth. Computing $\bar
{Q}_{b,n}^h$ for a given bandwidth is the only point in our simulations
where we do not treat $g_0$ as known. For a candidate blip function
estimate $\bar{Q}_b$, define the loss
\[
L_{\bar{Q}_0,g_0}(\bar{Q}_b) (o)\triangleq \biggl( \biggl[
\frac
{2a-1}{g_0(a|w)}\bigl(y-\bar{Q}_0(a,w)\bigr) +
\bar{Q}_0(1,w) - \bar {Q}_0(0,w) \biggr] -
\bar{Q}_b(w) \biggr)^2.
\]
To save computation time, we behave as though $\bar{Q}_0$ and $g_0$ are
known when using the above loss. We selected the bandwidth $H_n$ using
$10$-fold cross-validation with the above loss function to select from
the candidates $h=(0.01,0.02,\ldots,0.20)$. We also behave as though
$\bar
{Q}_0$ and $g_0$ are known when estimating each $\tilde{D}_{n,j}$, so
that the function $\tilde{D}_{n,j}$ only depends on $O_1,\ldots,O_{j-1}$
through the estimate of the optimal rule. This is mostly for
convenience, as it saves on computation time and our estimate of the
optimal rule $d_0^*$ will still not stabilize, that is, our optimal
value estimators will still encounter the irregularity at exceptional
laws. Note that $g_0$ is known in an RCT, and subtracting and adding
$\bar{Q}_0$ in the definition of the loss function will only serve to
stabilize the variance of our cross-validated risk estimate. In
practice, one could substitute an estimate of $\bar{Q}_0$ and expect
similar results. We update our estimates $d_{n,j}$ and $\tilde{\sigma
}_{0,n,j}$ using the method discussed in Section~\ref{sec:nuischunks} with
$S=\frac{n-\ell_n}{\ell_n}$.

To explore the sensitivity to the choice of $\ell_n$, we also
considered $(n,\ell_n)$ pairs $(1000,100)$ and $(4000,400)$, where
these pairs are only considered where explicitly noted. To explore the
sensitivity of our estimators to permutations of our data, we ran our
estimator twice on each Monte Carlo draw, with the indices of the
observations permuted so that the online estimator sees the data in a
different order.

\subsection{Simulation TTP-E: Two time point simulation} \label{sec:ttpsim}
The simulation used in this section was described in Section~8.1.2 of
\citet{vanderLaanLuedtke2014}, though here we modify the distribution
slightly so that the second time point treatment has no effect on the outcome.

\subsection*{Data}
The data is generated as follows:
\begin{eqnarray*}
&&L_1(0),L_2(0)\stackrel{\mathrm{i.i.d.}} {\sim}
\operatorname{Uniform}(-1,1),
\\
&&A(0)|L(0)\sim\operatorname{Bernoulli}(1/2),
\\
&&U_1,U_2|A(0),L(0)\stackrel{\mathrm{i.i.d.}} {\sim}
\operatorname{Uniform}(-1,1),
\\
&&L_1(1)|A(0),L(0),U_1,U_2\sim
U_1\times\bigl(1.25A(0)+0.25\bigr),
\\
&&L_2(1)|A(0),L(0),L_1(1),U_1,U_2
\sim U_2\times\bigl(1.25A(0)+0.25\bigr),
\\
&&A(1)|A(0),\bar{L}(1)\sim\operatorname{Bernoulli}(1/2),
\\
&&Y|\bar{A}(1),\bar{L}(1)\sim\operatorname{Bernoulli} \bigl(0.4 + 0.0345A(0) b
\bigl(L(0)\bigr) \bigr),
\end{eqnarray*}
where $b (L(0) )\triangleq-0.8 - 3 (\operatorname{sgn}(L_1(0))
+L_1(0) )-L_2(0)^2$. The treatment decision at time point 0 may
rely on $L(0)$, and the treatment at time point $1$ may rely on $L(0)$,
$A(0)$ and $L(1)$.

\subsection*{Estimation methods}
As in the previous simulation, we assume that the treatment mechanism
is known and supply the online estimator with correct estimates of the
conditional mean outcome so that $\tilde{D}_{n,j}$ is random only
through the estimate of $d_0^*$ (see Supplementary Appendix B for a
definition of $\tilde{D}_{n,j}$ in the two time point case). Given a
training sample $O_1,\ldots,O_j$, our estimator of $d_0^*$ corresponds to
using the full candidate library of weighted classification and
blip-function based estimators listed in Table~2 of \citet
{LuedtkevanderLaan2014}, with the weighted log loss used to determine
the convex combination of candidates. We update our estimate $d_{n,j}$
and $\tilde{\sigma}_{0,n,j}$ using the method described in
Section~\ref{sec:nuischunks} with $S=\frac{n-\ell_n}{\ell_n}$.

\subsection{Comparison with the $m$-out-of-$n$ bootstrap}
We compared our approach to the $m$-out-of-$n$ bootstrap for the value
of an estimated rule as presented by \citet{Chakrabortyetal2014}. By the
theoretical results in Section~\ref{sec:discdnconv}, it is reasonable
to expect that the optimal rule estimate will perform well and that one
can obtain inference for the optimal value using these same CIs. We ran
the $m$-out-of-$n$ bootstrap on D-E, C-NE, and C-E, with the same
sample sizes given in {Table~\ref{tab:nelln}}. We
drew $500$ bootstrap samples per Monte Carlo simulation, where we did
$500$ Monte Carlo simulations per setting due to the burdensome
computation time.

\begin{figure}

\includegraphics{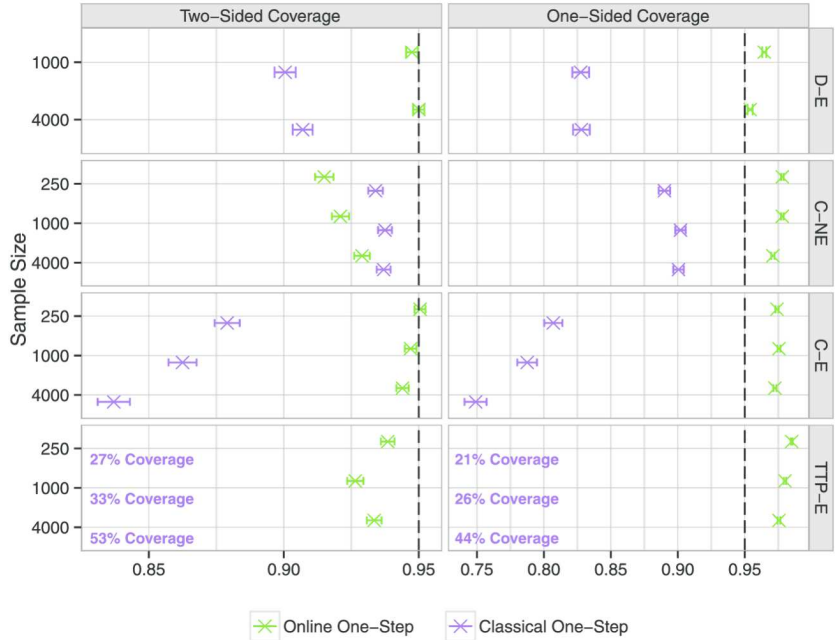}

\caption{Coverage of 95\% two-sided and one-sided (lower) CIs. The
online one-step estimator achieves near-nominal coverage for all of the
two-sided CIs and attains better than nominal coverage for the
one-sided CI. The classical one-step estimator only achieves
near-nominal coverage for C-NE. Error bars indicate 95\% CIs to account
for Monte Carlo uncertainty.}
\label{fig:covg}
\end{figure}

The $m$-out-of-$n$ bootstrap requires a choice of $m$, the size of each
nonparametric bootstrap sample. \citeauthor{Chakrabortyetal2014}
present a double bootstrap procedure for the two-time point case when
the optimal value is restricted to belong to a class of linear decision
functions. Because we do not restrict the set of possible regimes to
have linear decision functions, we instead set $m$ equal to $0.1n$,
$0.2n, \ldots, n$. When $n=1000$ and $m=100$, the NPMLE for D-E is
occasionally ill-defined due to empty strata. For these bootstrap
draws, we return the true optimal value, thereby (very slightly)
improving the coverage of the $m$-out-of-$n$ confidence intervals. We
will compare our procedure to the oracle regime, that is, the $m$ which
yields the shortest average CI length which achieves valid type I error
control. That is, we assume that one already knows the (on average)
optimal choice of $m$ from the outset.

\section{Simulation results} \label{sec:simresults}
\subsection{Online one-step compared to classical one-step}
Figure~\ref{fig:covg} shows the coverage attained by the online and
classical (non-online) one-step estimates of the optimal value. The
two-sided CIs resulting from the online estimator (nearly) attains
nominal coverage for all simulations considered, whereas the non-online
estimator only (nearly) attains nominal coverage for the
non-exceptional law in C-NE. The one-sided CIs from the online one-step
estimator attain proper coverage for all simulation settings. The
one-sided CIs from the non-online one-step estimates do not (nearly)
achieve nominal coverage in any of the simulations considered because
the rule is estimated on the same data as the optimal value. Thus,\vadjust{\goodbreak} we
expect to need a large sample size for the positive bias of the
non-online one-step to be negligible. In \citet{vanderLaanLuedtke2014},
we avoided this finite sample positive bias at non-exceptional laws by
using a cross-validated TMLE for the optimal value.

\begin{figure}

\includegraphics{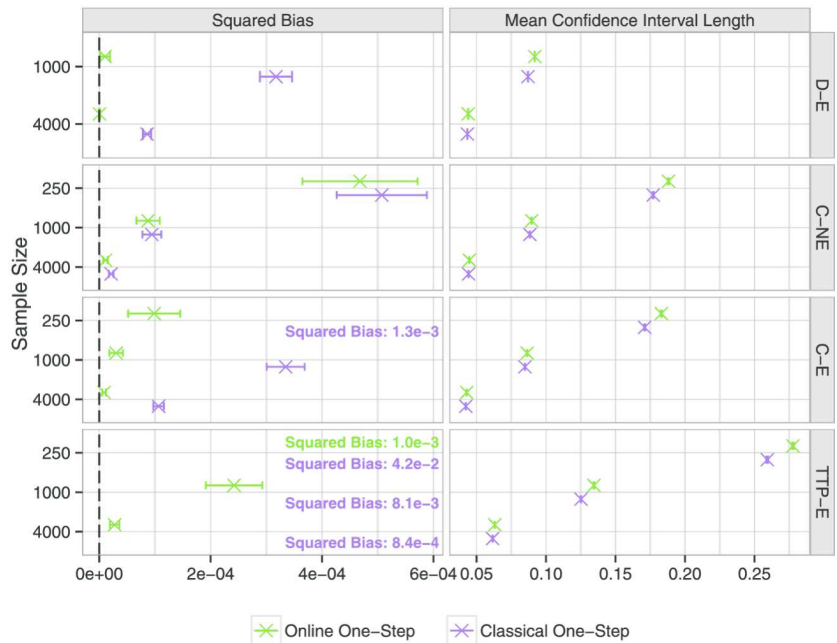}

\caption{Squared bias and 95\% two-sided CI lengths for the online and
classical one-step estimators, where the mean is taken across $2000$
Monte Carlo draws. The online estimator has lower squared bias than the
non-online estimator, while its mean CI length is only slightly longer.
Error bars indicate 95\% CIs to account for Monte Carlo uncertainty.}
\label{fig:bl}
\end{figure}

Figure~\ref{fig:bl} displays the squared bias and mean CI length across
the 2000 Monte Carlo draws. The online estimator consistently has lower
squared bias across all of our simulations. The online estimator was
negatively biased in all of our simulations, whereas the non-online
estimator was positively biased in all of our simulations. This is not
surprising: Theorem~\ref{thm:cilb} already implies that the online
estimator will generally be negatively biased in finite samples,
whereas the non-online estimator will generally be positively biased as
we have discussed.

\subsection{Online one-step compared to $m$-out-of-$n$ bootstrap}
\begin{figure}

\includegraphics{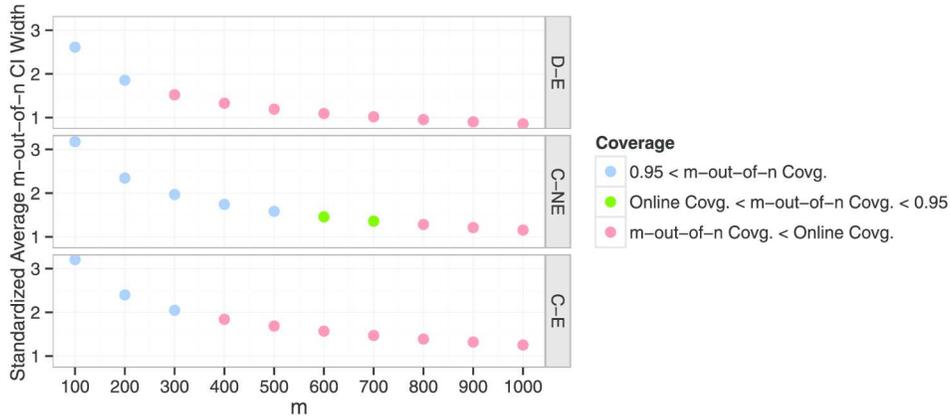}

\caption{Performance of the $m$-out-of-$n$ bootstrap at sample size
$1000$. The vertical axis shows the average CI width divided by the
average CI width of the online one-step CI. That is, any vertical axis
value above $1$ indicates the $m$-out-of-$n$ bootstrap has on average
wider CIs than the online one-step CI.}
\label{fig:mn}
\end{figure}

{Figure~\ref{fig:mn}} shows that our estimator
outperforms the $m$-out-of-$n$ bootstrap D-E and C-E for all choices of
$m$ considered at sample size $1000$. That is, our CI achieves
near-nominal coverage and is essentially always narrower than the CI
from the $m$-out-of-$n$, except when $m$ is very nearly equal to the
sample size. When $m$ is nearly equal to the sample size, coverage is
low for D-E and C-E: for $m=n$, the coverage is respectively 77\% and
65\%. When $m$ does achieve near-nominal coverage, the average CI width
is between 1.5 and 2 times larger than the average width of the online
one-step CIs. For C-NE, the estimation problem is regular and the
bootstrap performs (reasonably) well as expected by theory.
Nonetheless, so does the online procedure, and the online procedure
yields CIs of slightly shorter length for C-NE. The same general
results hold at other sample sizes, which we show in Figure A.3 in
Supplementary Appendix C.

One might argue that our oracle procedure is not truly optimal, since
in principle one could select a different choice of $m$ for each Monte
Carlo draw. While a valid criticism, we believe the overwhelming
evidence in favor of the online estimator presented in {Figure~\ref
{fig:mn}} should convince users that the online
approach will typically outperform any selection of $m$ at exceptional
laws. As $m$ is selected to be much less than $n$ at exceptional laws,
the $m$-out-of-$n$ will typically yield wider CIs than our procedure.
Given that our procedure has achieved near-nominal coverage at all
simulation settings, it seems hard to justify such a loss in precision.

\subsection{Sensitivity to permutations of the data and choice of
\texorpdfstring{$\ell_n$}{elln}}
While we would hope that our estimator is not overly sensitive to the
order of the data, the online estimator we have proposed necessarily
relies on a data ordering. Figure A.2 in Supplementary Appendix C
demonstrates how the optimal value estimates vary for C-E when the
estimator is computed on two permutations of the same data set. Our
point estimates are somewhat sensitive to the ordering of the data, but
this sensitivity decreases as sample size grows. We computed two-sided
CIs based on the two permuted data sets. We found that either both or
neither CI covered the true optimal value in 94\%, 94\% and 93\% of the
Monte Carlo draws at sample sizes 250, 1000 and 4000, respectively. For
C-NE, either both or neither CI covered the true optimal value in 91\%,
93\% and 95\% of the Monte Carlo draws at sample sizes 250, 1000 and
4000, respectively.

Different choices of $\ell_n$ did not greatly affect the coverage in
C-E and C-NE. Increasing $\ell_n$ for C-E decreased the coverage by
less than 1\% for sample sizes 1000 and 4000. Increasing $\ell_n$ for
C-NE increased the coverage by less than 1\% for sample sizes 1000 and
4000. Mean CI length increased predictably based on the increased value
of $\sqrt{n-\ell_n}$: for $n=1000$, increasing $\ell_n$ from $25$ to
$100$ increased the CI length by a multiplicative factor of $\sqrt
{\frac
{1000-25}{1000-100}}\approx1.04$. Similarly, increasing $\ell_n$ from
$100$ to $400$ increased the CI length by a factor of $1.04$ for $n=4000$.



\section{Discussion and future work} \label{sec:disc}
We have accomplished two tasks in this work. The first was to establish
conditions under which we would expect that regular root-$n$ rate
inference is possible for the mean outcome under the optimal rule. In
particular, we completely characterize the pathwise differentiability
of the optimal value parameter. This characterization on the whole
agrees with that implied by \citet{RobinsRotnitzky14}, but differs in a
minor fringe case where the conditional variance of the outcome given
covariates and treatment is zero. This fringe case may be relevant if
everyone in a strata of baseline covariates is immune to a disease
(regardless of treatment status) but are still included in the study
because experts are unaware of this immunity {a priori}. In
general, however, the two characterizations agree.

The remainder of our work shows that one can obtain an asymptotically
unbiased estimate of and a CI for the optimal value under reasonable
conditions. This estimator uses a slight modification of the online
one-step estimator presented by \citet{vanderLaanLendle2014}. Under
reasonable conditions, this estimator will be asymptotically efficient
among all RAL estimators of the optimal value at non-exceptional laws
in the nonparametric model where the class of candidate treatment
regimes is unrestricted. The main condition for the validity of our CI
is that the value of one's estimate of the optimal rule converges to
the optimal value at a faster than root-$n$ rate, which we show is
often a reasonable assumption. The lower bound in our CI is valid even
if this condition does not hold. 

We confirmed the validity of our approach using simulations. Our
two-sided CIs attained near-nominal coverage for all simulation
settings considered, while our lower CIs attained better than nominal
coverage (they were conservative) for all simulation settings
considered. Our CIs were of a comparable length to those attained by
the non-online one-step estimator. The non-online one-step estimator
only attained near-nominal coverage for the simulation which used a
non-exceptional data generating distribution, as would be predicted by theory.

In future work, we hope to mitigate the sensitivity of our estimator to
the order of the data. While we showed in our simulations that the
effect of permutations was minor, this property may be unappealing to
some. Such problems occur for many sample-split estimators; however,
one often has the option of estimating the parameter on several
permutations of the data and then averaging these estimates together.
The typical argument for averaging sample split estimates together is
that the estimator is asymptotically linear, that is, approximately an
average of a deterministic function applied to each of the $n$ i.i.d.
observations. Under mild conditions, we have an estimator which,
properly scaled, is equivalent to a sum of \textit{random} functions
applied to the $n$ observations, where these functions rely only on
past observations, making it impossible to apply this typical argument.
Further study is needed to determine if one can remove finite sample
noise from this estimator without affecting its asymptotic behavior.


Unsurprisingly, there is still more work to be done in estimating CIs
for the optimal rule. While we have shown that the lower bound from our
CI maintains nominal coverage under mild conditions, the upper bound
requires the additional assumption that the optimal rule is estimated
at a sufficiently fast rate. We observed in our simulations that the
non-online estimate of the optimal value had positive bias for all
settings. This is to be expected if the optimal rule is chosen to
maximize the estimated value, and can easily be explained analytically
under mild assumptions. It may be worth replacing the upper bound
${UB}_n$ in our CI by something like $\max\{{UB}_n,\psi_n(d_n)\}$,
where $\psi_n(d_n)$ is a non-online one-step estimate or TMLE of the
optimal value. One might expect that the upper bound $\psi_n(d_n)$ will
dominate the maximum precisely when the optimal rule is estimated
poorly. 


Finally, we note that our estimation strategy is not limited to
unrestricted classes of optimal rules. One could replace our
unrestricted class with, for example, a~parametric working model for
the blip function and expect similar results. This is because the
pathwise derivative of $P\mapsto E_{P_0}[Y_{d(P)}]$, which treats the
$P_0$ in the expectation subscript as known, will typically be zero
when $d(P)$ is an optimal rule in some class and does not fall on the
boundary of that class (with respect to some metric). Such a result
does not rely on $d(P)$ being a unique optimal rule. When the pathwise
derivative of $P\mapsto E_{P_0}[Y_{d(P)}]$ is zero, one can often prove
something like Theorem~\ref{thm:suffA5}, which shows that the value of the
estimated rule converges to the optimal value at a faster than root-$n$
rate under conditions.

Here, we considered the problem of developing a confidence interval for
the value of an unknown optimal treatment rule in a non-parametric
model. Under reasonable conditions, our proposed optimal value
estimator provides an interpretable and statistically valid approach to
gauging the effect of implementing the optimal individualized treatment
regime in the population.

\section*{Acknowledgments}
The authors would like to thank Sam Lendle for suggesting the
permutation analysis in our simulation, Robin Mejia and Antoine Chambaz
for greatly improving the readability of the document, and the
reviewers for helpful comments.

\begin{supplement}[id=suppA]
\stitle{Supplementary appendices: Proofs and extension to multiple
time point case}
\slink[doi]{10.1214/15-AOS1384SUPP} 
\sdatatype{.pdf}
\sfilename{aos1384\_supp.pdf}
\sdescription{Supplementary Appendix A contains all the proofs of all
of the results in the main text. Supplementary Appendix B contains an
outline of the extension to the multiple time point case. Supplementary
Appendix C contains additional figures referenced in the main text.}
\end{supplement}

%
%



\printaddresses
\end{document}